\documentclass[11pt, a4paper]{amsart}
\usepackage{amsmath,amsfonts,amssymb,amsthm}
\usepackage{graphicx} 
\usepackage{enumitem}

\title{On a generalization of a result of Kleitman}
\author{Ryan R. Martin}
\address{Iowa State University, Ames, Iowa, USA} 
\thanks{Research partially supported by Simons Foundation Collaboration Grant for Mathematicians \#709641}
\email{rymartin@iastate.edu}
\author{Bal\'azs Patk\'os}
\address{HUN-REN Alfr\'ed R\'enyi Institute of Mathematics} 
\thanks{Research supported by NKFIH under grant FK132060} 
\email{patkos@renyi.hu}
\date{}

\newtheorem{thm}{Theorem}
\newtheorem{cor}[thm]{Corollary}
\newtheorem{lem}[thm]{Lemma}
\newtheorem{prop}[thm]{Proposition}

\newtheorem{obs}[thm]{Observation}

\newcommand{\F}{{\mathcal F}}
\newcommand{\cE}{{\mathcal E}}
\newcommand{\cF}{{\mathcal F}}
\newcommand{\cG}{{\mathcal G}}
\newcommand{\cM}{{\mathcal M}}
\newcommand{\cH}{{\mathcal H}}
\newcommand{\cS}{{\mathcal S}}
\newcommand{\cA}{{\mathcal A}}
\newcommand{\cB}{{\mathcal B}}

\newcommand{\Kn}{\mathop{}\!\mathrm{Kn}}
\newcommand{\vex}{\mathop{}\!\mathrm{vex}}

\newcommand{\emb}{\mathop{}\!\mathrm{emb}}

\begin{document}

\begin{abstract}
A classical result of Kleitman determines the maximum number $f(n,s)$ of subsets in a family $\cF\subseteq 2^{[n]}$ of sets that do not contain distinct sets $F_1,F_2,\dots,F_s$ that are pairwise disjoint in the case $n\equiv 0,-1$ (mod $s$). 
Katona and Nagy determined the maximum size of a family of subsets of an $n$-element set that does not contain $A_1,A_2,\dots,A_t,B_1,B_2,\dots,B_t$ with $\bigcup_{i=1}^t A_i$ and $\bigcup_{i=1}^t B_i$ being disjoint. 
In this paper, we consider the problem of finding the maximum number $\vex(n,K_{s\times t})$ in a family $\cF\subseteq 2^{[n]}$ without sets $F^1_1,\dots,F^1_t,\dots,F^s_1,\dots,F^s_t$ such that $G_j=\bigcup_{i=1}^tF^j_i$ $j=1,2,\dots,s$ are pairwise disjoint. We determine the asymptotics of $2^n-\vex(n,K_{s\times t})$ if $n\equiv -1$ (mod $s$) for all $t$, and if $n\equiv 0$ (mod $s$), $t\ge 3$ and show that in this latter case the asymptotics of the $t=2$ subcase is different from both the $t=1$ and $t\ge 3$ subcases.
\end{abstract}

\maketitle

\section{Introduction}

An old and not completely settled problem in extremal set theory is the following: What is the maximum number of subsets in a family $\cF$ of an $n$-element set such that $\cF$ does not contain $s$ sets that are pairwise disjoint? 
The case $s=2$ is that of non-uniform intersecting families and the solution is an easy observation in the seminal paper of Erd\H os, Ko, and Rado \cite{EKR}. 
The case $s\ge 3$ and $n\equiv 0,-1$ $(\text{mod}\ s)$ was solved by Kleitman \cite{Kl} (see Theorem~\ref{Kthm} below). 
The remaining subcase of $s=3$ was first solved by Quinn in his thesis \cite{Q}, and later reproved and extended by Frankl and Kupavskii \cite{FK, FK2} for general $s$ and $n\equiv -2$ $(\text{mod}\ s)$. 

Answering a question of K\"orner (via private communication), Katona and Nagy \cite{KN} determined the maximum size of a family of subsets of an $n$-element set that does not contain two pairs of sets $A_1,A_2,B_1,B_2$ with $(A_1\cup A_2)$ and $(B_1\cup B_2)$ being disjoint. 
They also determined the asymptotics of the extremal number if the pairs are replaced by $t$-tuples $A_1,\dots,A_t,B_1\dots,B_t$ with $\bigl(\bigcup_{i=1}^tA_i\bigr)$ and $\bigl(\bigcup_{i=1}^tB_i\bigr)$ being disjoint.

There exists a general framework, introduced in \cite{GP}, to treat these problems. 
We will use the standard notation $[n]=\{1,2,\dots,n\}$, $2^X=\{S: S\subseteq X\}$, and $\binom{X}{k}=\{S\subseteq X: |S|=k\}$. 
The much studied \textit{Kneser graph} $\Kn(n,m)$ is the graph with vertex set $\binom{[n]}{m}$ such that $F$ and $F'$ are joined by an edge if and only if they are disjoint. 
A similar but much less studied graph, the \textit{Kneser cube}, is our main interest: $\Kn_n$ is the graph with vertex set $2^{[n]}$ such that $F$ and $F'$ are joined by an edge if and only if they are disjoint. 
For a graph $G$, we say that a family $\cF\subseteq 2^{[n]}$ is $G$-free if $\Kn_n[\cF]$ (the subgraph of the Kneser cube induced by $\cF$) is $G$-free. 

The above results are all special cases of the problem of determining $\vex(n,G)$, the size of largest $G$-free family $\cF\subseteq 2^{[n]}$ with $G$ being the complete graph $K_s$ or the complete bipartite graph $K_{t,t}$. 
Determining the value of $\vex(n,G)$ is an instance of vertex Tur\'an problems that have been well-studied, including in \cite{AT,AKSz,GetalA,Getal,JE,JT,Ko}.

\hskip 0.3truecm

Let us now formally state the result of Kleitman mentioned above. 
We will use the notation $\binom{n}{\le m}$, $\binom{n}{<m}$, $\binom{n}{\ge m}$, $\binom{n}{>m}$ for $\sum_{i=0}^m\binom{n}{i}$, $\sum_{i=0}^{m-1}\binom{n}{i}$, $\sum_{i=m}^n\binom{n}{i}$, $\sum_{i=m+1}^n\binom{n}{i}$, respectively.


\begin{thm}[Kleitman \cite{Kl}]\label{Kthm}
        ~\\ \vspace{-0.5cm}
        \begin{enumerate}[label={\it(\alph*)}]
       \item 
        if $n=sm-1$, then $2^n-\vex\bigl(n,K_s\bigr)=\binom{n}{\leq m-1}$, \label{Kthm:mod1}
        \item 
        if $n=sm$, then $2^n-\vex\bigl(n,K_s\bigr)=\binom{n}{m}-\binom{n-1}{m}+\binom{n}{\leq m-1}$. \label{Kthm:mod0}
    \end{enumerate}
   Moreover, the families that show that the value of $\vex\bigl(n,K_s\bigr)$ is sharp is shown by $\cF_1=\bigcup_{i=m}^n\binom{[n]}{i}$ for the case $n=sm-1$ and by $\cF_2=\binom{[n-1]}{m} \cup \bigcup_{i=m+1}^n\binom{[n]}{i}$ for the case $n=sm$.
\end{thm}

Observe that we express the solution in terms of $2^n-\vex\bigl(n,G\bigr)$. 
Hence, the goal is to minimize $2^{[n]}-|\F|$ over all $G$-free families $\F$. 
It should be noted that if the family $\F$ is $G$-free with $F\subsetneq F', F\in \cF, F'\notin \F$, then $\F\setminus \{F\}\cup \{F'\}$ is also $G$-free. 
Therefore when bounding the size of a largest $G$-free family, one can (and we will) assume that $\F$ is an upward-closed family. 
So $2^n-\vex(n,G)$, roughly speaking, tells how many lower layers of the Kneser cube should be dropped so that what remains is a $G$-free family.

We address the problem of determining $2^n-\vex(n,K_{s\times t})$ asymptotically in the cases $n\equiv -1,0\pmod{s}$ addressed by Kleitman's theorem, where $K_{s\times t}$ is the complete $s$-partite graph with all parts of size $t$. The case $s=2$, as mentioned above, was solved by Katona and Nagy~\cite{KN} and the case $t=1$ is simply the case of $K_s$, i.e. Theorem~\ref{Kthm} above. 
Our main result, Theorem~\ref{kst}, roughly states that if $s\ge 3$ and $t\geq 2$ then, except in the intriguing case of $n=sm$ and $t=2$, the bounds of Kleitman can be ``pushed down by one level''. 

\begin{thm}\label{kst}
Let $s\ge 3$.
\begin{enumerate}[label={\it(\alph*)}]
    \item 
    If $n=sm-1$, then for any $t\ge 2$, we have $2^n-\vex\bigl(n,K_{s\times t}\bigr)=\bigl(1+o(1)\bigr)\binom{n}{\le m-2}$. \label{kst:mod1}
    \item 
    If $n=sm$, then
    \begin{enumerate}[label={\it(\roman*)}]
        \item 
        for any $t\ge 3$, we have $2^n-\vex\bigl(n,K_{s\times t}\bigr)=\bigl(1+o(1)\bigr)\bigl[\binom{n}{m-1}-\binom{n-1}{m-1}+\binom{n}{\le m-2}\bigr]$, \label{kst:mod0:3}
        \item 
        for $t=2$, there exists $\delta_s>0$ such that
        \begin{align*} \binom{n}{\leq m-1}-\bigl(1+o(1)\bigr)\beta_s\binom{n}{m-1}\ge 2^n-\vex(n,K_{s\times 2})\ge \binom{n}{\leq m-1}-\Bigl(\frac{s-1}{s}-\delta_s+o(1)\Bigr)\binom{n}{m-1} , 
        \end{align*}
        where \label{kst:mod0:2}
    \begin{align*}
        \beta_s=\frac{s-1}{s-2}\biggl(1-\Bigl(1-\frac{s-2}{s^2-s}\Bigr)^{s-1}\biggr).
    \end{align*} 
    \end{enumerate} \label{kst:mod0}
\end{enumerate}
\end{thm}

Let us elaborate on how the bounds of Theorem \ref{kst} compare to those of Theorem \ref{Kthm}. 
The~\ref{Kthm:mod1} parts of these theorems are easier to summarize: Theorem \ref{Kthm}\ref{Kthm:mod1} states that the complement of a $K_s$-free family on an $sm-1$ sized ground set $X$ must contain at least as many sets as the number of subsets of $X$ with size at most $m-1$. 
Theorem \ref{kst}\ref{kst:mod1} pushes this result one layer lower: for any $t\ge 2$, the complement of a $K_{s\times t}$-free family on an $sm-1$ sized ground set $X$ must contain at least asymptotically as many sets as the number of subsets of $X$ with size at most $m-2$.

Something similar happens in the case when the ground set $X$ has size $sm$ and $t\geq 3$. 
In this case, Theorem \ref{Kthm}\ref{Kthm:mod0} states that the complement of a $K_s$-free family must contain at least as many sets as the number of subsets of $X$ of size at most $m-1$ and those $m$-subsets that contain a fixed element $x$ of $X$. 
Theorem \ref{kst}\ref{kst:mod0}\ref{kst:mod0:3} pushes this result by one layer asymptotically but only for $t\ge 3$. 
This is not because our methods are not strong enough to prove such a statement for the case $t=2$, but because in that case $2^n-\vex(n,K_{s\times 2})$ has different asymptotics than both $2^n-\vex(n,K_s)$ and $2^n-\vex(n,K_{s\times 3})$ as shown by Theorem \ref{kst}\ref{kst:mod0}\ref{kst:mod0:2}. 
The asymptotics of $2^{sm}-\vex(sm,K_{s\times 2})$ are still unknown, though.

\vskip 0.2truecm

To see how the results of $n=sm$ and $t=2$ differ from that of $n=sm$ and $t=1$, the upper bound of Theorem \ref{kst}\ref{kst:mod0}\ref{kst:mod0:2} shows that when changing $t$ from $1$ to $2$ one can add to $\F_2$  (the extremal family of the case $t=1$) not only the rest of the $m^{\rm th}$ layer, but also a positive proportion of the missing $(m-1)^{\rm st}$ layer of the cube.

On the other hand, the lower bound of Theorem \ref{kst}\ref{kst:mod0}\ref{kst:mod0:2} minus the expression in Theorem \ref{kst}\ref{kst:mod0}\ref{kst:mod0:3} gives
\begin{align*}
    \Bigl(2^n-\vex\bigl(n,K_{s\times 2}\bigr)\Bigr) - \Bigl(2^n-\vex\bigl(n,K_{s\times 3}\bigr)\Bigr) =\bigl(1+o(1)\bigr)\delta_s\binom{n}{m-1} .
\end{align*}

%
%

\vskip 0.2truecm

In Section \ref{sec:general}, we make some general observations on $\vex(n,G)$. By connecting $\vex(n,G)$ to forbidden subposet problems, Gerbner and the second author obtained the following result for bipartite graphs $G$.

\begin{thm}[\cite{GP}]
    For any bipartite graph $G$ containing at least one edge, we have $\bigl(1-o(1)\bigr)\binom{n}{\lfloor n/2 \rfloor }\ge \vex(n,G)-2^{n-1}\ge \bigl(1/2-o(1)\bigr)\binom{n}{\lfloor n/2 \rfloor }$.
\end{thm}

For a graph $G$ and positive integer $n$, we denote by $\emb(n,G)$ the maximum value of $m$ such that the Kneser graph $\Kn(n,m)$ contains a copy of $G$. 

 \begin{thm}\label{order}
     For any graph $G$ with chromatic number at least 3, we have $2^n-\vex(n,G)=\Theta_G\bigl(\binom{n}{\emb(n,G)}\bigr)$.
 \end{thm}

There exist several papers in the literature concerning subgraphs of Kneser graphs (see e.g. \cite{D,K,PT}), but homomorphisms to Kneser graphs are more often studied. The \textit{fractional chromatic number} $\chi_f(G)$ of a graph $G$ can be defined in several ways, one of which is $$\chi_f(G):=\inf\biggl\{\frac{a}{b}:\exists ~\text{a homomorphism}\ h:G\rightarrow \Kn(a,b)\biggr\}.$$
If $G$ is fixed and $n$ tends to infinity, then there is simple connection between embeddings and homomorphisms.

 \begin{prop}\label{triv}
     For any non-empty graph $G$, we have $\lim_{n\rightarrow \infty}\frac{\emb(n,G)}{n}=\frac{1}{\chi_f(G)}$.
 \end{prop}

 Putting together Theorem \ref{order} and Proposition \ref{triv}, one can obtain results how hard it is to determine the order of magnitude of $2^n-\vex(n,G)$. More precisely, if we are only interested in $\lim_{n\rightarrow\infty}\frac{\log[2^n-\vex(n,G)]}{n}$, that is the logarithm of the real number $c=c(G)$ with the property that for any $\varepsilon>0$ the inequalities $(c-\varepsilon)^n<2^n-\vex(n,G)<(c+\varepsilon)^n$ hold, then 
\begin{equation}\label{eq1}
   \log c(G)= \lim_{n\rightarrow\infty} \frac{\log [2^n-\vex(n,G)]}{n}=\lim_{n\rightarrow\infty} \frac{\log \binom{n}{\emb(n,G)}}{n}=H(1/\chi_f(G)),
\end{equation}
where $H(x)=-x\log x-(1-x)\log(1-x)$ is the binary entropy function.

Lund and Yannakis proved \cite{LY} that there exists a $\delta> 0$ such that there does not exist a polynomial time approximation algorithm of the fractional chromatic number that achieves the ratio $|V(G)|^\delta$ unless $P = NP$. Using (\ref{eq1}), one immediately obtains that approximating $c(G)$ within the ratio an appropriately chosen function $f(|V(G)|)$ is impossible unless $P=NP$.


We close this section with a basic inequality of binomial coefficients that will be useful later. 

\begin{prop}
    For any $m<n/2$, $\sum_{i=0}^m\binom{n}{i} \leq \frac{n-m}{n-2m}\binom{n}{m}$. \label{prop:geometric}
\end{prop}

\begin{proof}
    \begin{align*}
        \sum_{i=0}^m\binom{n}{i} 
        &= \sum_{i=0}^m\binom{n}{m}\frac{m!(n-m)!}{i!(n-i)!} = \sum_{j=0}^m\binom{n}{m}\frac{m!(n-m)!}{(m-j)!(n-m+j)!} \\
        &= \binom{n}{m}\sum_{j=0}^m\prod_{k=0}^{j-1}\frac{m-k}{n-m+j-k} \leq \binom{n}{m}\sum_{j=0}^m\Bigl(\frac{m}{n-m}\Bigr)^j \leq \binom{n}{m}\frac{1}{1-\frac{m}{n-m}} .  
    \end{align*}
\end{proof}

\subsection{Notation}

For a set $X\subseteq[n]$, we denote the \textit{shadow of $X$} with $\partial X=\{X\setminus x : x\in X\}$. For a family $\cF$ of sets, we write $\cF^i$ to denote the subfamily of all $i$-element sets of $\cF$.

\section{Constructions}

In this section, we describe the constructions showing the lower bounds in Theorem \ref{kst}.

\medskip

If $n=sm-1$, one can take all subsets of $[n]$ of size at least $m-1$. 
Any $t\ge 2$ distinct sets have a union of size at least $m$, so there is no space for $s$ pairwise disjoint such unions. 
This shows $\vex(sm-1,K_{s\times t})\ge \binom{sm-1}{\ge m-1}$. 
Note that if $t\ge 3$, then one can add a family $\cM$ of sets of size $m-2$ with the property that any $(m-1)$-set contains at most $t-2$ members of $\cM$. 
However, such an $\cM$ would have size at most $\frac{t-2}{m-1}\binom{n}{m-2}=o\bigl(\binom{n}{m-2}\bigr)$.

\medskip

If $n=sm$ and $t\ge 3$, then one can consider the family $\cF$ of all sets of size at least $m$ together with sets of size $m-1$ not containing a fixed element $x$ of $[n]$. 
The union of $t$ distinct sets from $\cF$ has size at least $m$, so if there existed $s$ such pairwise disjoint unions, call them $U_1,\ldots,U_s$, then each has size exactly $m$, their union is $[n]$, and thus $x$ must be contained in one such union, say $x\in U_1$. 
Since $|U_1|=m$, it can only contain one set of size at least $m$ (that is, $U_1$ itself) and thus must contain two other sets of size $m-1$.
However, since all the $(m-1)$-sets in $\cF$ fail to contain $x$, there is only one subset of $U_1$ of size $m-1$ in $\cF$; that is, $U_1\setminus x$, a contradiction. 
As a result, $\cF$ is $K_{s\times t}$-free for $t\ge 3$ and thus $\vex(sm,K_{s\times t})\ge \binom{sm}{\ge m}+\binom{sm-1}{m-1}$. 
Therefore, $2^n-\vex(sm,K_{s\times t})\le \binom{sm}{\le m-2}+\binom{sm-1}{m-2}$. 
Again, for $t\ge 4$, a lower order term family can be added.

\medskip

All that remains is the case $n=sm$ and $t=2$. We start with an observation.

\begin{obs}
    Let $\cF$ be an upward-closed family. 
    It has a copy of $K_{s\times 2}$ if and only if there are pairwise-disjoint sets $F_1,\ldots,F_s$ such that $\sum_{i=1}^s |F_i|\leq n-s$. \label{obs:s2}
\end{obs}

\begin{proof}
    Suppose first that $\cF$ contains $F_1,F_2,\dots,F_s$ that are pairwise disjoint and $\sum_{i=1}^s|F_i|\le n-s$. Then there exist $x_1,x_2,\dots,x_s\in [n]\setminus \bigcup_{i=1}^sF_i$. As a result, the sets $F_1, F_1\cup \{x_1\},F_2,F_2\cup \{x_2\},\dots,F_s,F_s\cup \{x_s\}$ form a copy of $K_{s\times 2}$. 
    In the other direction, suppose $F_{1,1},F_{1,2},F_{2,1},F_{2,2},\dots,F_{s,1},F_{s,2}$ form a copy of $K_{s\times 2}$ with $|F_{i,1}|\le |F_{i,2}|$ for all $1\le i\le s$ and the sets $F_{i,1}\cup F_{i,2}$ are pairwise disjoint. 
    Then, writing $F'_{i,2}=F_{i,1}\cup F_{i,2}$, $F_{1,1},F'_{1,2},F_{2,1},F'_{2,2},\dots,F_{s,1},F'_{s,2}$ also form a copy of $K_{s\times 2}$ with $F_{i,1}\subsetneq F'_{i,2}$. But the $F_{i,2}'$'s are pairwise disjoint with $\bigl|F_{i,1}\bigr|\leq \bigl|F_{i,2}'\bigr|-1$. Hence, $\sum_{i=1}^s\bigl|F_{i,1}\bigr|\leq \sum_{i=1}^s\bigl|F_{i,2}'\bigr|-s=n-s$.
\end{proof}

\begin{lem}\label{constr}
    For every integer $s\geq 3$, there exists an upward-closed $K_{s\times 2}$-free family in $2^{[sm]}$ of size $\binom{sm}{\geq m}+\bigl(1+o(1)\bigr)\beta_s\binom{sm}{m-1}$, where 
    \begin{align*}
        \beta_s=\frac{s-1}{s-2}\biggl(1-\Bigl(1-\frac{s-2}{s^2-s}\Bigr)^{s-1}\biggr).
    \end{align*}
    \label{lem:construction}
\end{lem}

\begin{proof}
    The construction is simple: For $i=1,\ldots,s-1$, let $\cF_i$ be the set of all sets of size $m-i$ containing at least $i$ elements of $[s-1]$. The family is $\cS:=\binom{[sm]}{\geq m}\cup\bigcup_{i=1}^{s-1}\cF_i$. 

Suppose $S_1,\ldots,S_s$ are pairwise disjoint. Then, as $\bigl|S_i\cap [s-1]\bigr|\ge \max\bigl\{m-\bigl|S_i\bigr|,0\bigr\}$, we have
    \begin{align*}
        sm - \sum_{i=1}^s \bigl|S_i\bigr| 
        = \sum_{i=1}^s\bigl(m-\bigl|S_i\bigr|\bigr) 
        &\le \sum_{i=1}^s \max\bigl\{m-\bigl|S_i\bigr|,0\bigr\} 
        \le \Biggl|\bigcup_{i=1}^s \bigl(S_i\cap [s-1]\bigr)\Biggr|
        \le s-1 \\
        sm - s + 1 
        &\le \sum_{i=1}^s \bigl|S_i\bigr| , 
    \end{align*}
    so by Observation~\ref{obs:s2}, $\cS$ is $K_{s\times 2}$-free.

    As to its size: 
    \begin{align*}
        \bigl|\cF_i\bigr|   
        &=  \sum_{j=i}^{s-1} \binom{s-1}{j} \binom{sm-s+1}{m-i-j} \\
        &=  \bigl(1+o(1)\bigr)\binom{sm}{m-1} \sum_{j=i}^{s-1} \binom{s-1}{j}\frac{(s-1)^{s-i-j}}{s^{s-1}} .
    \end{align*}

    Thus,
    \begin{align*}
        \beta_s     
        &=  \sum_{i=1}^{s-1} \sum_{j=i}^{s-1} \binom{s-1}{j} \frac{(s-1)^{s-i-j}}{s^{s-1}} \\
        &=  \frac{(s-1)^s}{s^{s-1}} \sum_{j=1}^{s-1} \sum_{i=1}^{j} \binom{s-1}{j} \Bigl(\frac{1}{s-1}\Bigr)^{i+j} \\
        &=  \frac{(s-1)^s}{s^{s-1}} \sum_{j=1}^{s-1} \binom{s-1}{j} \Bigl(\frac{1}{s-1}\Bigr)^{j} \cdot \frac{1}{s-2}\biggl(1-\Bigl(\frac{1}{s-1}\Bigr)^{j}\biggr) \\
        &=  \frac{(s-1)^s}{s^{s-1}(s-2)} \sum_{j=1}^{s-1} \binom{s-1}{j} \Bigl(\frac{1}{s-1}\Bigr)^{j} - \frac{(s-1)^s}{s^{s-1}(s-2)} \sum_{j=1}^{s-1} \binom{s-1}{j} \biggl(\frac{1}{(s-1)^2}\biggr)^{j} \\
        &=  \frac{(s-1)^s}{s^{s-1}(s-2)} \Biggl[\biggl(1+\frac{1}{s-1}\biggr)^{s-1} - 1\Biggr] - \frac{(s-1)^s}{s^{s-1}(s-2)} \Biggl[\biggl(1+\frac{1}{(s-1)^2}\biggr)^{s-1} - 1\Biggr] \\
        &=  \frac{(s-1)^s}{s^{s-1}(s-2)} \Biggl[\biggl(\frac{s}{s-1}\biggr)^{s-1} - \biggl(\frac{s^2-2s+2}{(s-1)^2}\biggr)^{s-1}\Biggr] ,
    \end{align*}
    which simplifies to the statement in the lemma. 
\end{proof}

Note that if $\cF'\subseteq \binom{[sm]}{m-1}$ is $K_s$-free, then $\cF'\cup \binom{[sm]}{\ge m}$ is $K_{s\times 2}$-free.
Indeed, in a copy of $K_{s\times 2}$ the union of each part would have size at least $m$, but as $\cF'$ is $K_s$-free, one part cannot contain sets from $\cF'$ and so the union of that part has size at least $m+1$, so altogether the parts have a total size at least $sm+1$ -- a contradiction. 
According to the celebrated Erd\H os Matching Conjecture \cite{E,Fra1.5,FK3}, in this particular case, the largest such $\cF$ should be $\cF_1'$, which is the construction of Lemma \ref{constr}. 
For this particular $\cF'$, the family $\cF'\cup \binom{[sm]}{\ge m}$ could still be extended while keeping it $K_{s\times 2}$-free. 
But if the Erd\H os Matching Conjecture fails, in theory it is possible that for the largest family $\cF'$ with the above property, $\cF'\cup \binom{[sm]}{\ge m}$ is maximally $K_{s\times 2}$-free. 
Let $\alpha_s$ be the expression forming the lower bound for the Erd\H{o}s Matching Conjecture; that is $\alpha_s=1-\bigl(1-\frac{1}{s}\bigr)^{s-1}$. See Table~\ref{table:construction} for a comparison between $\alpha_s$ and $\beta_s$.
\begin{table}
    \begin{tabular}{|r|r|r|} \hline
        $s$     &   $\alpha_s$  & $\beta_s$ \\ \hline\hline
        $3$     &   $5/9\approx 0.555556$  &   $11/18\approx 0.611111$ \\
        $4$     &   $37/64=0.578125$    &   $91/144\approx 0.631944$ \\
        $5$     &   $369/625=0.5904\;\;\;\,$    &   $25493/40000=0.637325$ \\ \hline
        $6$     &   $0.598122$    &   $0.638818$ \\
        $7$     &   $0.603431$    &   $0.639087$ \\
        $8$     &   $0.607304$    &   $0.638926$ \\ \hline
        $9$     &   $0.610256$    &   $0.638615$ \\
        $10$     &   $0.612580$    &   $0.638263$ \\
        $11$     &   $0.614457$    &   $0.637912$ \\ \hline
        $12$     &   $0.616005$    &   $0.637579$ \\
        $13$     &   $0.617303$    &   $0.637271$ \\
        $14$     &   $0.618408$    &   $0.636989$ \\ \hline
        $15$     &   $0.619360$    &   $0.636732$ \\
        $16$     &   $0.620188$    &   $0.636497$ \\
        $17$     &   $0.620915$    &   $0.636283$ \\ \hline
        $18$     &   $0.621558$    &   $0.636087$ \\
        $19$     &   $0.622132$    &   $0.635907$ \\
        $20$     &   $0.622646$    &   $0.635743$ \\ \hline
    \end{tabular}
    \caption{Table of coefficients $\beta_s$ for the construction in Lemma~\ref{lem:construction} compared to the lower bound for the Erd\H{o}s Matching Conjecture, $\alpha_s$.}
    \label{table:construction}
\end{table}

\section{Upper bounds on $K_{s\times t}$-free families}
 
In this section, we will prove the upper bounds of Theorem \ref{kst} on $\vex(n,K_{s\times t})$. First, we will prove the bounds corresponding to Theorem~\ref{kst}\ref{kst:mod1} and Theorem~\ref{kst}\ref{kst:mod0}\ref{kst:mod0:3}.
This will require only Kleitman's results (Theorem \ref{Kthm}) and some easy observations (Proposition \ref{good} below). 
Then to obtain the upper bound in Theorem~\ref{kst}\ref{kst:mod0}\ref{kst:mod0:2}, we will need some more tools and a closer look at the proof of Theorem \ref{Kthm}.

From here on, let $\cF\subseteq 2^{[n]}$ be a $K_{s\times t}$-free family with $n=sm$ or $n=sm-1$ and let $\cH:=2^{[n]}\setminus \cF$. 
If $F\in \cF,G\notin \cF$ with $F\subsetneq G$, then $\cF\setminus \{F\}\cup \{G\}$ is also $K_{s\times t}$-free, so to bound $\vex(n,K_{s\times t})$ we can and will assume that $\cF$ is upward closed and $\cH$ is downward closed.

Let $\pi$ be a partition of $[n]$ into $p\ge s$ parts. 
Then, the starting point of the proof of Kleitman's theorem is that the $K_{s\times 1}$-property of $\cF$ implies that at least $p-s+1$ parts of $\pi$ belong to $\cH$. 
This is not necessarily true for $t\ge 2$. Instead, we say that a set $X$ is \textit{good (with respect to $\cF$)} if $\bigl|\partial X\cap \cF\bigr|\le t-2$. 
The analog with respect to Theorem~\ref{Kthm} of the above observation is that the $K_{s\times t}$-property of $\cF$ implies that at least $p-s+1$ parts of $\pi$ are good. The family of good sets is denoted by $\cG$. 

\begin{obs}\label{obs:goodparts}
    Let $\cF$ be a $K_{s\times t}$-free family. We say that a set $X$ is \textrm{good (with respect to $\cF$)} if $\bigl|\partial X\cap \cF\bigr|\le t-2$. For any partition $\pi$ of $[n]$ into $p\geq s$ parts, at least $p-s+1$ parts of $\pi$ are good. The family of good sets is denoted by $\cG$. 
\end{obs}


As $\cH$ is downward closed,  $\cG\supseteq\cH$. Recall that $\cG^i$ are the sets in $\cG$ of size $i$ and $\cH^i$ is defined analogously. 

\begin{prop}With $\cG$ and $\cH$ as defined above:\label{good}
\begin{enumerate}[label={\it(\roman*)}]
    \item $\cG$ is downward closed. \label{it:downward}

    \item $\cE:=2^{[n]}\setminus \cG$ is $K_s$-free. \label{it:KsFree}

    \item $\bigl|\cG^i\bigr|(i-t+2)\le\bigl|\cH^{i-1}\bigr|(n-i+1)$ and $\displaystyle\frac{\bigl|\cG^i\bigr|}{\binom{n}{i}}\le \biggl(1+O\Bigl(\frac{1}{i}\Bigr)\biggr)\frac{\bigl|\cH^{i-1}\bigr|}{\binom{n}{i-1}}$. \label{it:GH}
\end{enumerate}
\end{prop}

\begin{proof}
    Suppose $x\in G\in \cG$. To show \ref{it:downward}, we need to prove that $G':=G\setminus \{x\}$ is good. 
    If $G'\in \cH$, then as $\cH$ is downward closed, all sets in $\partial G'$ are in $\cH\subseteq \cG$, and thus $|\partial G'\cap \cF|=0\le t-2$. 
    If $G'\notin \cH$, then for any $y\in G'$ we have $G'\setminus \{y\}\subset G\setminus \{y\}$. 
    As there exist at most $t-2$ such $y$ for which $G\setminus \{y\}\in \cF$ and $\cF$ is upward closed, there can be at most $t-2$ elements $y$ with $G'\setminus \{y\}\in \cF$ so $G'$ is good.

    To see \ref{it:KsFree}, $\cE$ consists of the non-good subsets, i.e. for every $E\in \cE$ there exist at least $t-1$ sets in $\partial E\cap \cF$. 
    But then a copy of $K_s$ in $\cE$ together with $t-1$ sets from the shadows (note that as the sets corresponding to $K_s$ are pairwise disjoint, therefore the shadows are pairwise disjoint) form a $K_{s\times t}$ in $\cF$, a contradiction.

    Finally for \ref{it:GH}, consider the bipartite graph $B$ with parts $\cG^i$ and $\cH^{i-1}$ with $G,H$ forming an edge in $B$ if and only if $H\subset G$. 
    The number of edges in $B$ is at least $\bigl|\cG^i\bigr|(i-t+2)$ by the definition of being a good set, and is at most $\bigl|\cH^{i-1}\bigr|(n-i+1)$ as every $(i-1)$-set has $n-i+1$ supersets of size $i$. 
    Thus, $\bigl|\cG^i\bigr|(i-t+2)\le\bigl|\cH^{i-1}\bigr|(n-i+1)$ and rearranging yields the inequality
    $$ \frac{\bigl|\cG^i\bigr|}{\binom{n}{i}}\le \Bigl(1+\frac{t-2}{i-t+2}\Bigr)\frac{\bigl|\cH^{i-1}\bigr|}{\binom{n}{i-1}}.$$
\end{proof}

The next theorem, together with the constructions from the previous section, establishes Theorem \ref{kst}\ref{kst:mod1} and Theorem~\ref{kst}\ref{kst:mod0}\ref{kst:mod0:3}.

\begin{thm}
    Suppose that $\cF\subseteq 2^{[n]}$ is a $K_{s\times t}$-free family with $s\ge 3$, $t\ge 2$. Then
    \begin{itemize}
        \item 
        if $n=sm$, then $\displaystyle |\cF|\le \Bigl(1+O\Bigl(\frac{1}{m}\Bigr)\Bigr) \biggl[\frac{s-1}{s}\binom{n}{m-1}+\sum_{i=m}^n\binom{n}{i}\biggr]$ and this bound is asymptotically sharp if $t\ge 3$,
        \item 
        if $n=sm-1$, then $\displaystyle |\cF|\le \Bigl(1+O\Bigl(\frac{1}{m}\Bigr)\Bigr) \sum_{i=m}^n\binom{n}{i}$ and this is asymptotically sharp.
    \end{itemize}
\end{thm}

\begin{proof}
Let $\cF\subseteq 2^{[n]}$ be as in the statement of the theorem, define $\cG$ to be the family of good sets and $\cH=2^{[n]}\setminus \cF, \cE=2^{[n]}\setminus \cG$. 
By Proposition \ref{good}\ref{it:KsFree}, $\cE$ is $K_s$-free and thus by Theorem \ref{Kthm}, we obtain
\begin{itemize}
    \item if $n=sm$, then
    $\displaystyle |\cG|\ge \frac{1}{s}\binom{n}{m}+\binom{n}{\le m-1}$ and thus $\displaystyle \sum_{i=m-m^{2/3}}^n|\cG^i|\ge \frac{1}{s}\binom{n}{m}+\sum_{i=m-m^{2/3}}^{m-1}\binom{n}{i}$;
    \item if $n=sm-1$, then
    $\displaystyle |\cG|\ge \binom{n}{\le m}$ and thus $\displaystyle \sum_{i=m-m^{2/3}}^n|\cG^i|\ge \sum_{i=m-m^{2/3}}^{m}\binom{n}{i}$.
    \end{itemize}

Applying Proposition \ref{good}\ref{it:GH} to all $\cG^i$ with $i\ge m-m^{2/3}$ we obtain that
\begin{itemize}
    \item 
    if $n=sm$, then    
    \begin{align*}
        |\cH| &\ge  \sum_{i=m-m^{2/3}-1}^n|\cH^i|\ge \frac{m-m^{2/3}-t+2}{sm-m+m^{2/3}+1}\sum_{i=m-m^{2/3}}^n|\cG^i| \\
        &\ge  \biggl(\frac{1}{s-1}-O\Bigl(\frac{1}{m^{1/3}}\Bigr)\biggr) \biggl[\frac{1}{s}\binom{n}{m}+\sum_{i=m-m^{2/3}}^{m-1}\binom{n}{i}\biggr] \\
        &\ge  \biggl(1-O\Bigl(\frac{1}{m^{1/3}}\Bigr)\biggr) \biggl[\frac{1}{s}\binom{n}{m-1}+\sum_{i=m-m^{2/3}}^{m-2}\binom{n}{i}\biggr],
    \end{align*}
    \item 
    if $n=sm-1$, then
    \begin{align*}
        |\cH| &\ge  \sum_{i=m-m^{2/3}-1}^n|\cH^i|  \ge  \frac{m-m^{2/3}-t+2}{sm-m+m^{2/3}+1}\sum_{i=m-m^{2/3}}^n|\cG^i| \\
        &\ge  \biggl(\frac{1}{s-1}-O\Bigl(\frac{1}{m^{1/3}}\Bigr)\biggr) \sum_{i=m-m^{2/3}}^{m-1}\binom{n}{i} \\
        &\ge  \biggl(1-O\Bigl(\frac{1}{m^{1/3}}\Bigr)\biggr) \sum_{i=m-m^{2/3}}^{m-2}\binom{n}{i} .
    \end{align*}
\end{itemize}
Taking complements give the required bounds on $\cF$ noting that $\sum_{i=0}^{m-m^{2/3}}\binom{n}{i}=o\bigl(\frac{1}{m}\binom{n}{m}\bigr)$ both for $n=sm$ and $n=sm-1$.
\end{proof}

It remains to prove the upper bound of Theorem \ref{kst}\ref{kst:mod0}\ref{kst:mod0:2}. To do this, we need to redo its proof, apply Proposition \ref{good} and add an extra tool.

Kleitman established some inequalities for complements of upward closed $K_s$-free families which we list below in Lemma~\ref{lem:Kleitman}. 
By Proposition \ref{good}~\ref{it:downward} and~\ref{it:KsFree}, it is the case that each of those inequalities are valid for $\cG$. 
Using Proposition \ref{good}~\ref{it:GH}, we will turn them into inequalities about $\cH$. 
To be able to state the inequalities in Lemma~\ref{lem:Kleitman}, we need to introduce the following: $P_r$ denotes the set of those partitions of $[n]$ that contain exactly $r$ good parts. If $\pi$ is an ordered partition with sizes $i_1,i_2,\dots,i_p$ (we say that such a partition has type $\pi$), then the number of ordered partitions of type $\pi$ is $n(\pi)=\frac{n!}{\prod_{j=1}^pi_j!}$. The ratio of all ordered partitions of type $\pi$ that contain exactly $r$ good sets is denoted by $\rho_r(\pi)$. Observation~\ref{obs:goodparts} implies 
$$\rho_j(\pi)=0 ~\text{for}\ j<p-s+1 ~\text{and}\ \sum_{r=0}^p\rho_r(\pi)=\sum_{r=p-s+1}^p\rho_r(\pi)=1. $$

\begin{lem}[Kleitman \cite{Kl}]\label{lem:Kleitman}
    Let $\pi_e$ denote an equipartition of $[sm]$ into $s$ parts. If $\cA\subseteq 2^{[sm]}$ is an upward closed $K_s$-free family, then for $\cB=2^{[sm]}\setminus \cA$ the following inequalities hold:
    \begin{enumerate}[label={\it(\roman*)}]
        \item 
        $\displaystyle s\frac{\bigl|\cB^m\bigr|}{\binom{sm}{m}}=\sum_{j=1}^sj\rho_j\bigl(\pi_e\bigr)=1+\sum_{j=1}^s(j-1)\rho_j\bigl(\pi_e\bigr)$, 
        \item 
        for $j\ge 1$: $$\frac{\bigl|\cB^{m-j}\bigr|}{\binom{sm}{m-j}}+(s-1)\frac{|\cB^{m+1}|}{\binom{sm}{m+1}}\ge \sum_{r=j+1}^s\frac{r}{j}\rho_r\bigl(\pi_e\bigr)+\sum_{r=1}^j\rho_j(\pi_e)=1-\sum_{r=j+1}^s\Bigl(1-\frac{r}{s}\Bigr)\rho_r\bigl(\pi_e\bigr),$$       
        \item 
        for $j\ge 0$: $$\frac{\bigl|\cB^{m-s+1-j}\bigr|}{\binom{sm}{m-s+1-j}}+(s-1)\frac{\bigl|\cB^{m+1}\bigr|}{\binom{sm}{m+1}}\ge 1.$$ 
    \end{enumerate}
\end{lem}

Applying Lemma~\ref{lem:Kleitman} with $\cB=\cG$, using Proposition \ref{good}~\ref{it:GH} and rearranging, we obtain the following corollary.

\begin{cor}\label{cor}
     Let $\pi_e$ denote an equipartition of $[sm]$ into $s$ parts. As long as $j=o(m)$, the following inequalities hold:
    \begin{enumerate}[label={\it(\roman*)}]
        \item 
        $\displaystyle \bigl|\cH^{m-1}\bigr|\ge \biggl(1 - O\Bigl(\frac{1}{m}\Bigr)\biggr)\frac{1}{s}\Biggl[\binom{sm}{m-1}+\binom{sm}{m-1}\sum_{j=1}^s(j-1)\rho_j(\pi_e)\Biggr]$, \label{it:cor:gen}
        \item 
        for $j\ge 1$: $$\bigl|\cH^{m-1-j}\bigr| + \frac{(s-1)\bigl|\cH^{m}\bigr|}{\binom{sm}{m}}\binom{sm}{m-1-j} \ge \biggl(1-O\Bigl(\frac{1}{m}\Bigr)\biggr) \biggl[\binom{sm}{m-1-j}-\binom{sm}{m-1-j}\sum_{r=j+1}^s\biggl(1-\frac{r}{s}\biggr)\rho_r\bigl(\pi_e\bigr)\biggr],$$ \label{it:cor:one} 
        \item 
        for $j\ge 0$: $$\bigl|\cH^{m-s-j}\bigr| + \frac{(s-1)\bigl|\cH^{m}\bigr|}{\binom{sm}{m}}\binom{sm}{m-s-j} \ge \biggl(1-O\Bigl(\frac{1}{m}\Bigr)\biggr) \binom{sm}{m-s-j}.$$ \label{it:cor:zero}
        \end{enumerate}
\end{cor}

We will need the following recent result of the authors.

\begin{thm}{\cite{MP}}\label{mp}
    For any $s\ge 3$ there exists $\varepsilon_s>0$ such that any $K_s$-free family $\cF\subseteq \binom{[n]}{h}$ with $n>sh$ has size at most $(1-\frac{1}{s}-\varepsilon_s)\binom{n}{h}$.
\end{thm}

Note that Observation \ref{obs:s2} implies that $\cF^{m-1}$ is $K_s$-free, and thus by Theorem \ref{mp}, we have $|\cH^{m-1}|\ge (\frac{1}{s}+\varepsilon_s)\binom{sm}{m-1}$. For any $0<\gamma <1$ one can take the linear combination of Corollary \ref{cor} \ref{it:cor:gen} (with coefficient $1-\gamma$) and this last inequality (with coefficient $\gamma$) to obtain

\begin{equation}\label{ineq}
\bigl|\cH^{m-1}\bigr|\ge \biggl(1 - O\Bigl(\frac{1}{m}\Bigr)\biggr)\frac{1}{s}\Biggl[(1+\delta)\binom{sm}{m-1}+(1-\gamma)\binom{sm}{m-1}\sum_{j=1}^s(j-1)\rho_j(\pi_e)\Biggr],    
\end{equation} 

where $\delta=s\varepsilon_s\gamma$.

\vskip 0.3truecm

Using the inequalities of Corollary \ref{cor} and (\ref{ineq}), we bound $\sum_{j=m'}^m\bigl|\cH^j\bigr|$ from below as follows.

\begin{align*}
    \sum_{j=m'}^m\bigl|\cH^j\bigr| 
    &= \bigl|\cH^m\bigr| + \bigl|\cH^{m-1}\bigr| + \sum_{j=m'}^{m-2}\bigl|\cH^j\bigr| \\
    &= \bigl|\cH^m\bigr| + \bigl|\cH^{m-1}\bigr| + \sum_{j=1}^{s-2}\bigl|\cH^{m-1-j}\bigr| + \sum_{j=0}^{m-s-m'}\bigl|\cH^{m-s-j}\bigr|\\
    &\ge \bigl|\cH^m\bigr| + \biggl(1 - O\Bigl(\frac{1}{m}\Bigr)\biggr)\frac{1}{s}\Biggl[(1+\delta)\binom{sm}{m-1}+(1-\gamma)\binom{sm}{m-1}\sum_{j=1}^s(j-1)\rho_j(\pi_e)\Biggr]\\
    &\hspace{12pt} -\frac{(s-1)\bigl|\cH^{m}\bigr|}{\binom{sm}{m}}\sum_{j=1}^{s-2}\binom{sm}{m-1-j} \\
    &\hspace{12pt} + \biggl(1-O\Bigl(\frac{1}{m}\Bigr)\biggr)\sum_{j=1}^{s-2} \Biggl[\binom{sm}{m-1-j}-\binom{sm}{m-1-j}\sum_{r=j+1}^s\biggl(1-\frac{r}{s}\biggr)\rho_r\bigl(\pi_e\bigr)\Biggr] \\
    &\hspace{12pt} -\frac{(s-1)\bigl|\cH^{m}\bigr|}{\binom{sm}{m}}\sum_{j=0}^{m-s-m'}\binom{sm}{m-s-j} + \biggl(1-O\Bigl(\frac{1}{m}\Bigr)\biggr)\sum_{j=0}^{m-s-m'}\binom{sm}{m-s-j} 
\end{align*}

We now rearrange terms:
\begin{align}
    \sum_{j=m'}^m\bigl|\cH^j\bigr| 
    &\ge \bigl|\cH^m\bigr| \Biggl[1-\frac{s-1}{\binom{sm}{m}}\sum_{j=m'}^{m-2}\binom{sm}{j}\Biggr] + \biggl(1 - O\Bigl(\frac{1}{m}\Bigr)\biggr) \Biggl[\frac{1+\delta}{s}\binom{sm}{m-1}+\sum_{j=m'}^{m-2}\binom{sm}{j}\Biggr] \nonumber \\
    &\hspace{12pt} + \biggl(1 - O\Bigl(\frac{1}{m}\Bigr)\biggr)\Biggl[\frac{1-\gamma}{s}\binom{sm}{m-1}\sum_{j=1}^s(j-1)\rho_j(\pi_e)\Biggr] \nonumber \\
    &\hspace{12pt} - \biggl(1-O\Bigl(\frac{1}{m}\Bigr)\biggr)\sum_{j=1}^{s-2}\biggl(\binom{sm}{m-1-j}\sum_{r=j+1}^s\biggl(1-\frac{r}{s}\biggr)\rho_r\bigl(\pi_e\bigr)\biggr) \nonumber \\
    &= \bigl|\cH^m\bigr|\Biggl[1-\frac{s-1}{\binom{sm}{m}}\sum_{j=m'}^{m-2}\binom{sm}{j}\Biggr] + \biggl(1 - O\Bigl(\frac{1}{m}\Bigr)\biggr)\Biggl[\frac{1+\delta}{s}\binom{sm}{m-1}+\sum_{j=m'}^{m-2}\binom{sm}{j}\Biggr] \label{eq:main} \\
    &\hspace{12pt} + \biggl(1 - O\Bigl(\frac{1}{m}\Bigr)\biggr)\Biggl[\sum_{r=2}^s\Biggl(\frac{(r-1)(1-\gamma)}{s}\binom{sm}{m-1}-\sum_{h=1}^{r-1}\biggl(1-\frac{r}{s}\biggr)\binom{sm}{m-1-h}\Biggr)\rho_r\bigl(\pi_e\bigr)\Biggr] \nonumber
\end{align}

Finally, we establish that the coefficients of $\bigl|\cH^m\bigr|$ and $\rho_r\bigl(\pi_e\bigr)$ are nonnegative. 
Using Proposition~\ref{prop:geometric} and $s\ge 3$,
\begin{align*}
    1-\frac{s-1}{\binom{sm}{m}}\sum_{j=m'}^{m-2}\binom{sm}{j}
    &\ge 1-\frac{s-1}{\binom{sm}{m}}\frac{sm-(m-2)}{sm-2(m-2)}\binom{sm}{m-2} = 1-\frac{(s-1)m(m-1)}{(sm-2m+4)(sm-m+1)} \\
    &> 1-\frac{1}{s-2} \geq 0 .
\end{align*}

For the coefficient of $\rho_r\bigl(\pi_e\bigr)$, we first consider the case $r=2$.
\begin{align*}
    \frac{(r-1)(1-\gamma)}{s}\binom{sm}{m-1}-\sum_{h=1}^{r-1}\biggl(1-\frac{r}{s}\biggr)\binom{sm}{m-1-h}
    &= \frac{1-\gamma}{s}\binom{sm}{m-1}-\frac{s-2}{s}\binom{sm}{m-2} \\
    &> \binom{sm}{m-1}\biggl(\frac{1-\gamma}{s}-\frac{s-2}{s}\cdot\frac{1}{s}\biggr)\ge 0,
\end{align*}
if $\gamma\le 2/s$.

\vskip 0.2truecm

For $r\ge 3$ using Proposition~\ref{prop:geometric} and the fact that $s\geq 3$, 
\begin{align*}
    \frac{(r-1)(1-\gamma)}{s}\binom{sm}{m-1}-\sum_{h=1}^{r-1}\biggl(1-\frac{r}{s}\biggr)\binom{sm}{m-1-h}
    &= \frac{(r-1)(1-\gamma)}{s}\binom{sm}{m-1}-\sum_{j=m-r}^{m-2}\biggl(1-\frac{r}{s}\biggr)\binom{sm}{j} \\
    &\ge \frac{(r-1)(1-\gamma)}{s}\binom{sm}{m-1}-\biggl(1-\frac{r}{s}\biggr)\frac{sm-m+2}{sm-2m+4}\binom{sm}{m-2} \\
    &= \binom{sm}{m-1}\Biggl[\frac{(r-1)(1-\gamma)}{s}-\biggl(1-\frac{r}{s}\biggr)\frac{m-1}{sm-2m+4}\Biggr] \\
    &> \binom{sm}{m-1}\Biggl[\frac{(r-1)(1-\gamma)}{s}-\biggl(1-\frac{r}{s}\biggr)\frac{1}{s-2}\Biggr].
\end{align*}
For $\gamma<(s-1)/(s-2)$, the expression in the brackets is monotone increasing in $r$, and for $r=3$ and $\gamma<(s-1)/(2s-4)$, this expression is non-negative as claimed.

\vskip 0.2truecm

Returning to \eqref{eq:main}, as their coefficients are non-negative as long as $\gamma$ is sufficiently small relative to $s$, we can lower bound  right hand side of \eqref{eq:main} by removing the expressions involving $|\cH^m|$ and $\rho_r\bigl(\pi_e\bigr)$. Letting $m'=m-\alpha+1$, where $\alpha=\alpha(m)=o(m)$ and using Proposition~\ref{prop:geometric} again, we obtain
\begin{align*}
    \sum_{j=m'}^m\bigl|\cH^j\bigr| 
    &\ge \biggl(1 - O\Bigl(\frac{1}{m}\Bigr)\biggr) \Biggl[\frac{1+\delta}{s}\binom{sm}{m-1}+\sum_{j=m'}^{m-2}\binom{sm}{j}\Biggr] \\
    &= \biggl(1 - O\Bigl(\frac{1}{m}\Bigr)\biggr) \Biggl[\frac{1+\delta}{s}\binom{sm}{m-1}+\sum_{j=0}^{m-2}\binom{sm}{j}-\sum_{j=0}^{m-\alpha}\binom{sm}{j}\Biggr] \\
    &\ge \biggl(1 - O\Bigl(\frac{1}{m}\Bigr)\biggr) \Biggl[\frac{1+\delta}{s}\binom{sm}{m-1}+\sum_{j=0}^{m-2}\binom{sm}{j}-\frac{sm-m+\alpha}{sm-2m+2\alpha}\binom{sm}{m-\alpha}\Biggr]
\end{align*}

It suffices to show that $\frac{sm-m+\alpha}{sm-2m+2\alpha}\binom{sm}{m-\alpha} = O\bigl(1/m\bigr)\binom{sm}{m-1}$. To that end,
\begin{align*}
    \frac{sm-m+\alpha}{sm-2m+2\alpha}\frac{\binom{sm}{m-\alpha}}{\binom{sm}{m-1}} 
    &< \frac{s-1}{s-2}\frac{\bigl(m-1\bigr)_{\alpha-1}}{\bigl(sm-m+\alpha\bigr)_{\alpha-1}} < \frac{s-1}{s-2}\cdot \frac{1}{(s-1)^{\alpha-1}},
\end{align*}
so any choice of $\alpha$ satisfying $\Omega(\log m)=\alpha=o(m)$ will do.

We note that since it suffices to choose any $\gamma\leq 2/s$, the value of $\delta_s$ stated in the theorem can be as large as $2\epsilon_s$.

\section{General results}
\label{sec:general}

In this section, we prove our results on $\vex(n,G)$ for general $G$.

\begin{proof}[Proof of Theorem \ref{order}]
Before starting the proof, let us note that by definition, for any $m\geq n/2$, the graph $\Kn(n,m)$ is empty and $\Kn(2m,m)$ is a matching. So if $G$ is a matching, then its fractional chromatic number is 2, so the statement of the theorem holds, otherwise we know that $\emb(n,G)<n/2$ so Proposition \ref{prop:geometric} applies with $m=\emb(n,G)$.
    
Observe  that $\emb(n,G)$ is linear in $n$ as $\Kn\bigl(n,\bigl\lfloor\frac{n}{|G|}\bigr\rfloor\bigr)$ contains a clique of size $|G|$ and thus $G$. 
Then, setting $\cF_{n,G}=\bigl\{F\subseteq [n] : |F|>\emb(n,G)\bigr\}$ we claim that $\cF_{n,G}$ is $G$-free, for $n$ large enough. 
Indeed, if $F_1,F_2,\dots F_{|G|}$ form a copy of $G$ in $\cF_{n,G}$, then one can pick distinct sets $F'_i\subseteq F_i$ of size $\emb(n,G)+1$ for all $1\le i\le |G|$ (here we use that $\emb(n,G)$ grows with $n$) contradicting the definition of $\emb(n,G)$. 
This shows $2^n-\vex(n,G)\le 2^n-\bigl|\cF_{n,G}\bigr|=\sum_{i=0}^{\emb(n,G)}\binom{n}{i}\le C\binom{n}{\emb(n,G)}$, where for the last inequality we used Proposition \ref{prop:geometric}.

To see a lower bound on $2^n-\vex(n,G)$, we use a standard double counting argument. By the symmetry of $K=\Kn(n,\emb(n,G))$, every set $F$ of size $\emb(n,G)$ is contained in the same $N$ number of copies of $G$. So the number of copies of $G$ in $K$ is $\frac{N\cdot \binom{n}{\emb(n,G)}}{|V(G)|}$. Therefore to avoid all copies of $G$, we need to remove at least $\frac{\binom{n}{\emb(n,G)}}{|V(G)|}$ vertices (that is, $\emb(n,G)$-sets) from $K$.
\end{proof}

\begin{proof}[Proof of Proposition \ref{triv}]
One of the required inequalities follows from the fact that every embedding is a homomorphism and so $\chi_f(G)\cdot\emb(n,G)\geq n$. 
To see the other inequality observe that if $h:G\rightarrow \Kn(a,b)$ is a homomorphism, then for any $k\ge 1$ one can define another homomorphism $h_k:G\rightarrow \Kn(ak,bk)$  by $i\in h(v)$ if and only if $\bigl[(i-1)k+1,ik\bigr]\subset h_k(v)$. 
Then if $h'_k:G\rightarrow \Kn(ak,bk-1)$ is defined by $h'_k(v)\in \partial h_k(v)$, then $h'_k$ is still a homomorphism, and by the marriage theorem $h'_k$ can be defined to be injective and thus an embedding.    
\end{proof}

\end{document}